\documentclass{amsart}
\usepackage{amsmath,amscd,xypic,amssymb,combelow,color,enumitem,graphicx,float,comment,mathtools,mathrsfs,tikz-cd,mathdots,amsbsy}

\newtheorem{theorem}[subsection]{Theorem}

\newtheorem{definition}[subsection]{Definition}
\newtheorem{claim}[subsection]{Claim}

\newtheorem{remark}[subsection]{Remark}

\def\fg{{\mathfrak{g}}}

\def\fsl{{\mathfrak{sl}}}
\def\fgl{{\mathfrak{gl}}}

\def\BA{{\mathbb{A}}}
\def\BC{{\mathbb{C}}}
\def\BK{{\mathbb{K}}}

\def\BP{{\mathbb{P}}}

\def\BQ{{\mathbb{Q}}}
\def\BZ{{\mathbb{Z}}}

\def\CO{{\mathcal{O}}}

\def\CS{{\mathcal{S}}}

\def\CV{{\mathcal{V}}}

\def\sym{\textrm{sym}}

\def\bn{{\boldsymbol{n}}}

\def\b0{{\boldsymbol{0}}}

\def\loccit{\emph{loc.~cit.~}}

\def\tua{\widetilde{U}_q^+(L\fg_A)}
\def\tuwa{\widetilde{U}_q^+(L\fg_{\widehat{A}})}
\def\tub{\widetilde{U}_q^+(L\fg_B)}
\def\tuc{\widetilde{U}_q^+(L\fg_C)}
\def\tud{\widetilde{U}_q^+(L\fg_D)}

\def\ua{U_q^+(L\fg_A)}
\def\uwa{U_q^+(L\fg_{\widehat{A}})}
\def\ub{U_q^+(L\fg_B)}
\def\uc{U_q^+(L\fg_C)}
\def\ud{U_q^+(L\fg_D)}

\def\oCS{\mathring{\CS}}
\def\fosp{\mathfrak{osp}}

\begin{document}

\title[QUIVERS WITH POTENTIAL FOR QUANTUM AFFINE SUPERALGEBRAS]{\Large{\textbf{Quivers with potential for quantum affine superalgebras}}} 

\author[Andrei Negu\cb t]{Andrei Negu\cb t}

\address{École Polytechnique Fédérale de Lausanne (EPFL), Lausanne, Switzerland \newline \text{ } \ \ Simion Stoilow Institute of Mathematics (IMAR), Bucharest, Romania} 

\email{andrei.negut@gmail.com}
\def\uppercasenonmath#1{}

\maketitle

\begin{abstract} We construct a surjective algebra homomorphism from the $K$-theoretic Hall algebras of certain quivers with potential to the (half) quantum affine and toroidal orthosymplectic superalgebras.  \end{abstract}

\section{Introduction}
\label{sec:intro}

\medskip 

A major topic in modern representation theory is the correspondence of algebraic structures that feature in different areas of mathematics and physics. In the present paper, we will be interested in comparing the following algebras that arise in quantum algebra, algebraic geometry, and mathematical physics (respectively):

\medskip 

\begin{itemize}[leftmargin=*]
	
\item quantum loop (super)algebras associated to quivers with parameters;

\medskip 

\item $K$-theoretic Hall algebras ($K$-HAs) of quivers with potential (\cite{P}, following \cite{KS});

\medskip 

\item BPS algebras arising from string theory on toric Calabi-Yau threefolds (\cite{HM}).

\end{itemize}

\medskip 

\noindent Correspondences between the three algebras above were studied in the mathematical physics literature: \cite{GLY, GY, NW} following the Yangian case studied in \cite{LY} (see also \cite{Bao} for the cases relevant to the present paper, and \cite{BY,L} for more general cases). Certain comparison theorems were also established mathematically in \cite{N Red, RSYZ, RSYZ 2, VV}. In the present paper, we discuss correspondences between the first two bullets in the list above, in the case of finite and affine orthosymplectic Lie superalgebras
\begin{equation}
\label{eqn:osp intro}
\fg = \fosp_{m|n} \quad \text{or} \quad \widehat{\fosp}_{m|n}
\end{equation}
subject to the non-special hypothesis of Definition~\ref{def:special} in the affine case \footnote{The case $\fg \in \fsl_{m|n}$ was treated in \cite{T}, while $\fg = \widehat{\fsl}_{m|n}$ for $m \neq n$ was treated in \cite{N Red}.}. For any such $\fg$, we compare the ``Drinfeld new" half of the quantum affinization
\begin{equation}
\label{eqn:quantum loop intro}
U_q^+(L\fg)
\end{equation}
defined in \cite{BFK,XZ} (finite type) and \cite{NT} (affine type) with the $T$-equivariant localized $K$-theoretic Hall algebra ($K$-HA for short; see Subsections~\ref{sub:general} and \ref{sub:mf} for details)
\begin{equation}
\label{eqn:k-ha intro}
K^T_{\text{crit}}(Q,W)_{\text{loc}}
\end{equation} 
of a certain quiver $Q$ with potential $W$. These quivers were defined for almost all types $\fosp_{m|n}$ in \cite{Bao}, though we will show that in the particular case that we call ``Type $C$ Case 3", one needs to choose the specific potential in formula \eqref{eqn:scary potential}. We also define suitable quivers with potential $(Q,W)$ for non-special (cf. Definition~\ref{def:special}) affine $\widehat{\fosp}_{m|n}$ types, and with this in mind, the following is our main result.

\medskip 

\begin{theorem}
\label{thm:main}

For any $\fg$ as in \eqref{eqn:osp intro} (subject to the non-special assumption of Definition~\ref{def:special} in the affine case), there is a surjective algebra homomorphism
\begin{equation}
\label{eqn:hom main}
K^T_{\emph{crit}}(Q,W)_{\emph{loc}} \twoheadrightarrow U_q^+(L\fg)
\end{equation}
where the quiver with potential $(Q,W)$ is defined in Subsection~\ref{sub:quiver with potential}.

\end{theorem}

\medskip 

\noindent Our main comparison between the two sides of \eqref{eqn:hom main} is via the shuffle algebra
\begin{equation}
\label{eqn:shuffle intro}
\CS^+ \subset \CV^+ = \bigoplus_{(n_i \geq 0)_{i \in I}} \text{Frac}(\text{Rep}_T)[z_{i1}^{\pm 1}, \dots, z_{in_i}^{\pm 1}]_{i \in I}^{\text{sym}}
\end{equation}
which is a particular set of Laurent polynomials endowed with the Feigin-Odesskii shuffle product (\cite{FO}). It was shown in \cite{NT} that the shuffle algebra $\CS^+$ is isomorphic to the codomain of \eqref{eqn:hom main}, so the main result of the present paper is to compare the shuffle algebra with the domain of \eqref{eqn:hom main}. There is an algebra homomorphism
\begin{equation}
\label{eqn:map to shuffle}
K^T_{\text{crit}}(Q,W)_{\text{loc}} \rightarrow \CV^+
\end{equation}
due to \cite{P}, inspired by \cite{KS}. Using an argument pioneered in \cite{Z}, we show that the image of the homomorphism above is equal to $\CS^+$, which implies Theorem~\ref{thm:main}. 

\medskip 
 
\subsection{Acknowledgements}  I would like to thank Shivang Jindal for numerous very interesting conversations about quivers with potential. I gratefully acknowledge the support of the Swiss National Science Foundation grant 10005316.

\bigskip 

\section{Quantum affine and toroidal superalgebras}
\label{sec:algebra}

\medskip 

\subsection{Decorated quivers}
\label{sub:quivers}

Consider a quiver (i.e. oriented graph) $Q$ with vertex set $I$ and arrow set $\Omega$. We call $Q$ \emph{decorated} if its arrows $\alpha$ are endowed with parameters $t_\alpha$ in a certain characteristic zero field $\BK$:
\begin{equation}
\label{eqn:graph}
\alpha : \begin{tikzcd}
i & j
\arrow["t_\alpha", from=1-1, to=1-2]
\end{tikzcd}
\end{equation}
Given such a decorated quiver, we may associate a collection of \emph{zeta functions}
\begin{equation}
\label{eqn:zeta}
\left\{ \zeta_{ij}(x) \sim \frac {\prod_{\alpha : i \mapsto j} (1-xt_{\alpha})}{(1-x)^{\delta_{ij}}} \right\}_{i,j\in I}
\end{equation}
where $\sim$ means that the two sides of the equation are equal up to multiplication by a monomial in $x$. The value of this monomial will not be important for the purposes of the present paper (in particular for the proof of Theorem~\ref{thm:main}), although it is important for matching the algebras considered in the present paper with similarly-named objects in the literature, and for constructing their Drinfeld doubles. We refer to \cite{NT} for the appropriate choices of monomials in the situations \eqref{eqn:osp intro}, and note that in the case of odd isotropic vertices of Lie superalgebras, these monomials may involve half-integer powers of $x$ (see \loccit for the cosmetic modifications one needs to make to the subsequent framework in order to deal with such powers).

\medskip 

\subsection{Algebras}
\label{sub:algebras}

\noindent One associates the \emph{(half) pre-quantum loop algebra} 
\begin{equation}
\label{eqn:pre}
\widetilde{\mathbf{U}}^+ = \BK \langle e_{i,d} \rangle_{i \in I, d\in \BZ} \Big/ \left(e_i(x)e_j(y) \zeta_{ji}\left(\frac yx\right) = e_j(y) e_i(x) \zeta_{ij}\left(\frac xy\right) \right)_{i,j \in I}
\end{equation}
to the zeta functions \eqref{eqn:zeta}, where 
\[
e_i(x) = \sum_{d \in \BZ} \frac {e_{i,d}}{x^d}
\]
and one interprets the quadratic relation in \eqref{eqn:pre} by clearing the denominators of the form $(x-y)^{\delta_{ij}}$ and then equating all monomials in $x,y$ in the two sides of the equation.  By analogy with \cite{E, FO}, we may also consider the \emph{big shuffle algebra} \footnote{As explained in \cite[Section 3.4]{NT}, the right-hand side of \eqref{eqn:big shuf} should in fact be redefined as
\begin{equation}
\label{eqn:big shuf correct}
\CV^+ = \bigoplus_{\bn = (n_i \geq 0)_{i \in I}} \underline{z}^{\bn} \BK[z_{i1}^{\pm 1}, \dots, z_{in_i}^{\pm 1}]^{\sym}_{i \in I}
\end{equation}
where $\underline{z}^{\bn} = \prod_{i \text{ odd isotropic}} \prod_{u=1}^{n_i} z_{iu}^{\frac {n_i-1}2}$ must be used due to various half-integral powers that are concealed in the notation $\sim$ of \eqref{eqn:zeta}. Since we are choosing to gloss over these powers in the present paper, we also choose to gloss over the discrepancy between \eqref{eqn:big shuf correct} and \eqref{eqn:big shuf}.}	
\begin{equation}
\label{eqn:big shuf}
\CV^+ = \bigoplus_{\bn = (n_i \geq 0)_{i \in I}} \BK\left[z_{i1}^{\pm 1}, \dots, z_{in_i}^{\pm 1}\right]^{\sym}_{i \in I}
\end{equation}
with the superscript ``sym" referring to color-symmetric Laurent polynomials, i.e. those which are symmetric in $z_{i1},\dots,z_{in_i}$ for each individual $i \in I$. Following \emph{loc. cit.}, we make \eqref{eqn:big shuf} into an associative algebra via the \emph{shuffle product} 
\begin{equation}
\label{eqn:shuffle product}
E(z_{i1},\dots,z_{in_i})_{i \in I} * E'(z_{i1},\dots,z_{in'_i})_{i \in I} = 
\end{equation}
$$
\text{Sym} \left[ \frac {E(z_{i1},\dots,z_{in_i})_{i \in I} E'(z_{i,n_i+1},\dots,z_{i,n_i+n_i'})_{i \in I}}{\bn!\bn'!} \prod_{i,j \in I} \prod_{a=1}^{n_i} \prod_{b=n_j+1}^{n_j+n_j'} \zeta_{ij} \left(\frac {z_{ia}}{z_{jb}} \right) \right] 
$$
where ``Sym" denotes summation over the $(\bn+\bn')! = \prod_{i \in I} (n_i+n_i')!$ permutations of the variables of each color taken separately. The relation between the algebras \eqref{eqn:pre} and \eqref{eqn:big shuf} is that there exists an algebra homomorphism
\begin{equation}
\label{eqn:tilde upsilon}
\widetilde{\Upsilon}^+ \colon \widetilde{\mathbf{U}}^+ \rightarrow \CV^+
\end{equation}
given by $e_{i,d} \mapsto z_{i1}^d$ for all $i \in I$ and $d \in \BZ$. The \emph{(half) quantum loop algebra} is 
\begin{equation}
\label{eqn:quad}
\mathbf{U}^+ = \widetilde{\mathbf{U}}^+ \Big / \text{Ker } \widetilde{\Upsilon}^+
\end{equation}
and the \emph{small shuffle algebra} is defined as
\begin{equation}
\label{eqn:shuf}
\oCS^+ = \text{Im } \widetilde{\Upsilon}^+
\end{equation}
By the first isomorphism theorem, we have an induced isomorphism
\begin{equation}
\label{eqn:upsilon}
\Upsilon^+ \colon \mathbf{U}^+ \xrightarrow{\sim} \oCS^+
\end{equation} 
A \emph{shuffle realization} of $\mathbf{U}^+$ (see \cite[Section 2]{NT} for an overview of the general framework) refers to the problem of finding an explicit subalgebra
\begin{equation}
\label{eqn:shuffle realization}
\CS^+ \subset \CV^+
\end{equation}
such that $\CS^+ = \oCS^+$. Historically, many such subalgebras were determined by \emph{wheel conditions}, as we will now explain. Suppose there exists a cycle $\bigcirc$ in the quiver $Q$
\begin{equation}
\label{eqn:wheel}
\begin{tikzcd}
&& {i_{k-1}} &&& {i_{k}} && \\ \\
\vdots &&&&&&& {i_{1}} \\ \\
&& {i_{3}} &&& {i_{2}}
\arrow["{\alpha_k}", from=1-3, to=1-6]
\arrow["{\alpha_1}", from=1-6, to=3-8]
\arrow[from=3-1, to=1-3]
\arrow["{\alpha_2}", from=3-8, to=5-6]
\arrow[from=5-3, to=3-1]
\arrow["{\alpha_3}", from=5-6, to=5-3]
\end{tikzcd}
\end{equation}
(vertex and arrow repetitions are allowed) such that
\begin{equation}
\label{eqn:wheel parameter}
t_{\alpha_1} \dots t_{\alpha_k} = 1
\end{equation}
but $t_{\alpha_{s+1}} \dots t_{\alpha_{s'}} \neq 1$ for all $1 \leq s < s' \leq k$ such that $i_s = i_{s'}$. 

\medskip 

\begin{definition}
\label{def:wheel}

We say that $E(z_{iu})_{i \in I}^{u \geq 1} \in \CV^+$ \emph{vanishes at the wheel} $\bigcirc$ of \eqref{eqn:wheel} if 
\begin{equation}
\label{eqn:satisfies wheel}
E \Big|_{z_{i_ku_k} = z_{i_{k-1}u_{k-1}} t_{\alpha_k}, \dots, z_{i_2u_2} = z_{i_1u_1} t_{\alpha_2}, z_{i_1u_1} = z_{i_ku_k} t_{\alpha_1}} = 0
\end{equation}
for all $u_1,\dots,u_k \geq 1$ such that $u_s \neq u_{s'}$ if $s \neq s'$ and $i_s = i_{s'}$. 

\end{definition}

\medskip 

\noindent It is elementary to show that the subset of $\CV^+$ consisting of elements which vanish at any collection of wheels $\bigcirc$ is a subalgebra with respect to the shuffle product \eqref{eqn:shuffle product}. Many shuffle realizations \eqref{eqn:shuffle realization} consist of polynomials which vanish at a certain collection of wheels $\bigcirc$ as in \eqref{eqn:wheel} \footnote{Though as we shall see in the statement of Theorem~\ref{thm:b,c,d 2}, some cases need a slightly stronger version of the vanishing condition \eqref{eqn:satisfies wheel}.}.

\medskip 

\subsection{The special linear superalgebras}
\label{sub:sec a}

We henceforth let $\BK = \BQ(q)$. Let us apply the construction in the previous Subsection to the quantum affinization of the special linear Lie superalgebra
\[
\fg_A = \fsl_{m|n}
\]
To do so, we must choose a collection of simple roots of $\fg_A$, which is well-known to be in one-to-one correspondence with so-called \emph{parity sequences}
\begin{equation}
\label{eqn:sequence}
(s_1,\dots,s_N) \in \{+1,-1\}^N
\end{equation}
where $N = m+n$, with $1$ occurring $m$ times and $-1$ occurring $n$ times. Let $\bar{I} = \{1,\dots,N-1\}$. Given a sequence \eqref{eqn:sequence}, consider the \emph{decorated Dynkin diagram}
\begin{equation}
\label{eqn:dynkin a}
\begin{tikzcd}
1 && 2 && \cdots && {N-2} && {N-1}
\arrow["q^{-s_{N-1}}", shift left=2pt, from=1-9, to=1-7]
\arrow["q^{-s_{N-1}}", shift left=2pt, from=1-7, to=1-9]
\arrow["q^{-s_{N-2}}", shift left=2pt, from=1-7, to=1-5]
\arrow["q^{-s_{N-2}}", shift left=2pt, from=1-5, to=1-7]
\arrow["q^{-s_3}", shift left=2pt, from=1-5, to=1-3]
\arrow["q^{-s_3}", shift left=2pt, from=1-3, to=1-5]
\arrow["q^{-s_2}", shift left=2pt, from=1-3, to=1-1]
\arrow["q^{-s_2}", shift left=2pt, from=1-1, to=1-3]
\arrow["q^{2s_{N-1}}"', dashed, from=1-9, to=1-9, looseness=4]
\arrow["q^{2s_{N-2}}"', dashed, from=1-7, to=1-7, looseness=4]
\arrow["q^{2s_2}"', dashed, from=1-3, to=1-3, looseness=4.5]
\arrow["q^{2s_1}"', dashed, from=1-1, to=1-1, looseness=4.5]
\end{tikzcd}
\end{equation}
where the convention is that the loop at vertex $i$ is present if and only if $s_i = s_{i+1}$ (in which case it is decorated by the parameter $q^{2s_i}$), and absent otherwise. By \eqref{eqn:zeta}, the quiver above determines the type $A$ zeta function
\begin{equation} 
\label{eqn:zeta a}
\zeta^{A}_{ij}(x) \sim
\begin{cases} 
\frac {1-xq^{2s_i}}{1-x} &\text{if } i = j \text{ has a loop} \\ 
\frac {1}{1-x} &\text{if } i = j \text{ does not have a loop} \\ 
1-xq^{-s_i} &\text{if } i = j+1 \\ 
1-xq^{-s_j} &\text{if } j = i+1 \\ 
1 &\text{otherwise} 
\end{cases} 
\end{equation}
for all $i,j \in \bar{I}$. Therefore, the (half) pre-quantum affine algebra is
\begin{equation}
\label{eqn:quadratic a}
\tua = \BK \langle e_{i,d} \rangle_{i \in \bar{I}, d\in \BZ} \Big/ \left(e_i(x)e_j(y) \zeta^{A}_{ji}\left(\frac yx\right) = e_j(y) e_i(x) \zeta^{A}_{ij}\left(\frac xy\right) \right)_{i,j \in \bar{I}}
\end{equation}
The corresponding shuffle algebra was studied in numerous papers, for example \cite{FJM, T}. The motivation for the definitions above comes from the following results.

\medskip 

\begin{theorem}
\label{thm:a 1}

(\cite{NT}) The kernel of the homomorphism $\widetilde{\Upsilon}^+$ is generated by the cubic and quartic Drinfeld-Serre relations defined in \cite{Y2}, and so
\[
\ua = \tua \Big/ \emph{Ker } \widetilde{\Upsilon}^+
\]
is the (Drinfeld new) half of the special linear quantum affine superalgebra. 

\end{theorem}

\medskip 

\begin{theorem}
\label{thm:a 2}

(\cite{T}) If we let $\CV^+_A$ and $\oCS^+_A$ be the constructions \eqref{eqn:big shuf} and \eqref{eqn:shuf} (respectively) for the zeta functions \eqref{eqn:zeta a}, then we have
\[
\oCS^+_A = \CS^+_A
\]
where $\CS^+_A \subset \CV^+_A$ consists of all Laurent polynomials which vanish at the wheels
\begin{equation}
\label{eqn:wheel finite even}
\begin{tikzcd}
i-1 && i &&&& i \\ \\
i &&&& i && i+1
\arrow["{q^{-s_{i}}}"', from=1-1, to=3-1]
\arrow["{q^{-s_i}}"', from=1-3, to=1-1]
\arrow["{q^{-s_{i}}}", from=1-7, to=3-7]
\arrow["{q^{2s_i}}"', from=3-1, to=1-3]
\arrow["{q^{2s_i}}", from=3-5, to=1-7]
\arrow["{q^{-s_i}}", from=3-7, to=3-5]
\end{tikzcd}
\end{equation}
for any $i \in \bar{I}$ such that $s_i = s_{i+1}$, and
\begin{equation}
\label{eqn:wheel finite odd}
\begin{tikzcd}
i-1 && i \\ \\
i && i+1
\arrow["{q^{-s_i}}", from=1-1, to=1-3]
\arrow["{q^{-s_{i+1}}}", from=1-3, to=3-3]
\arrow["{q^{-s_i}}", from=3-1, to=1-1]
\arrow["{q^{-s_{i+1}}}", from=3-3, to=3-1]
\end{tikzcd}
\end{equation}
for any $i \in \bar{I}$ such that $s_i \neq s_{i+1}$.

\end{theorem}

\medskip 

\subsection{The orthosymplectic superalgebras}
\label{sub:osp}

We will now go through the orthosymplectic Lie superalgebras, and classify them in types $B,C,D$ following \cite{NT}. The logic behind this classification has to do with the shape of the corresponding decorated Dynkin diagrams: in all the pictures below, the subdiagram with vertices $1,\dots,N-1-\delta_{\text{type}=D}$ is called the \emph{body} and is of the form \eqref{eqn:dynkin a} for a certain parity sequence \eqref{eqn:sequence}, while the subdiagram with vertices $N-1-\delta_{\text{type}=D},\dots,N$ is called the \emph{head} and resembles that of usual type $B,C,D$ Lie algebras. All the signs $\pm$ in a given diagram are one and the same. 

\medskip 

\noindent \textbf{Type $B$} recovers $\fosp_{2m+1|2n}$. There are two decorated Dynkin diagrams:

\medskip 

\begin{itemize}[leftmargin=*]
	
\item 

Case $1$: 
\begin{equation}
\label{eqn:dynkin b 1}
\begin{tikzcd}
\cdots && {N-2} && {N-1} && N
\arrow["q^{\pm 1}", shift left=2pt, from=1-7, to=1-5]
\arrow["q^{\pm 1}", shift left=2pt, from=1-5, to=1-7]
\arrow[shift left=2pt, from=1-5, to=1-3]
\arrow[shift left=2pt, from=1-3, to=1-5]
\arrow[shift left=2pt, from=1-3, to=1-1]
\arrow[shift left=2pt, from=1-1, to=1-3]
\arrow["q^{\mp 1}"', from=1-7, to=1-7, looseness=4]
\arrow[dashed, from=1-5, to=1-5, looseness=4]
\arrow[dashed, from=1-3, to=1-3, looseness=4]
\end{tikzcd}
\end{equation}

\medskip

\item 
Case $2$: 
\begin{equation}
\label{eqn:dynkin b 2}
\begin{tikzcd}
\cdots && {N-2} && {N-1} && N 
\arrow["q^{\pm 1}", shift left=2pt, from=1-7, to=1-5]
\arrow["q^{\pm 1}", shift left=2pt, from=1-5, to=1-7]
\arrow[shift left=2pt, from=1-5, to=1-3]
\arrow[shift left=2pt, from=1-3, to=1-5]
\arrow[shift left=2pt, from=1-3, to=1-1]
\arrow[shift left=2pt, from=1-1, to=1-3]
\arrow["q^{\mp 2}"', from=1-7, to=1-7, out=45, in=135, distance=10mm]
\arrow["q^{\pm 1}", from=1-7, to=1-7, out=315, in=225, distance=10mm]
\arrow[dashed, from=1-5, to=1-5, looseness=4]
\arrow[dashed, from=1-3, to=1-3, looseness=4]
\end{tikzcd}
\end{equation}

\end{itemize}

\medskip 

\noindent \textbf{Type $C$} recovers $\fosp_{2m|2n}$ that end with the longer even simple root, see~\cite[\S 3.1, diagrams ($C_N$)]{Y1}. There are three decorated Dynkin diagrams:

\medskip 

\begin{itemize}[leftmargin=*]

\item Case $1$: 
\begin{equation}
\label{eqn:dynkin c 1}
\begin{tikzcd}
\cdots &&  {N-2} && {N-1} && N
\arrow["q^{\pm 2}", shift left=2pt, from=1-7, to=1-5]
\arrow["q^{\pm 2}", shift left=2pt, from=1-5, to=1-7]
\arrow["q^{\pm 1}", shift left=2pt, from=1-5, to=1-3] 
\arrow["q^{\pm 1}", shift left=2pt, from=1-3, to=1-5] 
\arrow[shift left=2pt, from=1-3, to=1-1]
\arrow[shift left=2pt, from=1-1, to=1-3]
\arrow["q^{\mp 4}"', from=1-7, to=1-7, looseness=4]
\arrow["q^{\mp 2}"', from=1-5, to=1-5, looseness=4]
\arrow[dashed, from=1-3, to=1-3, looseness=4]
\end{tikzcd}
\end{equation}

\medskip 

\item Case $2$: 
\begin{equation}
\label{eqn:dynkin c 2}
\begin{tikzcd}
\cdots && N-3 && {N-2} && {N-1} && N
\arrow["q^{\pm 2}", shift left=2pt, from=1-9, to=1-7]
\arrow["q^{\pm 2}", shift left=2pt, from=1-7, to=1-9]
\arrow["q^{\mp 1}", shift left=2pt, from=1-7, to=1-5]
\arrow["q^{\mp 1}", shift left=2pt, from=1-5, to=1-7]
\arrow["q^{\pm 1}", shift left=2pt, from=1-5, to=1-3]
\arrow["q^{\pm 1}", shift left=2pt, from=1-3, to=1-5]
\arrow[shift left=2pt, from=1-3, to=1-1]
\arrow[shift left=2pt, from=1-1, to=1-3]
\arrow["q^{\mp 4}"', from=1-9, to=1-9, looseness=4]
\arrow[dashed, from=1-3, to=1-3, looseness=4]
\end{tikzcd}
\end{equation}

\medskip 

\item Case $3$: 
\begin{equation}
\label{eqn:dynkin c 3}
\begin{tikzcd}
\cdots && N-3 && {N-2} && {N-1} && N
\arrow["q^{\pm 2}", shift left=2pt, from=1-9, to=1-7]
\arrow["q^{\pm 2}", shift left=2pt, from=1-7, to=1-9]
\arrow["q^{\mp 1}", shift left=2pt, from=1-7, to=1-5]
\arrow["q^{\mp 1}", shift left=2pt, from=1-5, to=1-7]
\arrow["q^{\mp 1}", shift left=2pt, from=1-5, to=1-3]
\arrow["q^{\mp 1}", shift left=2pt, from=1-3, to=1-5]
\arrow[shift left=2pt, from=1-3, to=1-1]
\arrow[shift left=2pt, from=1-1, to=1-3]
\arrow["q^{\mp 4}"', from=1-9, to=1-9, looseness=4]
\arrow["q^{\pm 2}"', from=1-5, to=1-5, looseness=4]
\arrow[dashed, from=1-3, to=1-3, looseness=4]
\end{tikzcd}
\end{equation}

\end{itemize}

\medskip 

\noindent \textbf{Type $D$} recovers $\fosp_{2m|2n}$ for all possible Dynkin diagrams that end with a fork as in~\cite[\S 3.1, diagrams ($D_N$)]{Y1}. There are two decorated Dynkin diagrams:

\medskip

\begin{itemize}[leftmargin=*]

\item Case $1$: 
\begin{equation}
\label{eqn:dynkin d 1}
\begin{tikzcd}
&&&& {N-1} \\
\cdots && {N-2} \\
&&&& N
\arrow[shift left=2pt, from=2-1, to=2-3]
\arrow[shift left=2pt, from=2-3, to=2-1]
\arrow["q^{\pm 1}", shift left=2pt, from=2-3, to=1-5]
\arrow["q^{\pm 1}", shift left=2pt, from=1-5, to=2-3]
\arrow["q^{\pm 1}", shift left=2pt, from=2-3, to=3-5]
\arrow["q^{\pm 1}", shift left=2pt, from=3-5, to=2-3]
\arrow["q^{\mp 2}", swap, from=1-5, to=1-5, looseness=4]
\arrow["q^{\mp 2}", swap, from=3-5, to=3-5, looseness=4]
\arrow[dashed, from=2-3, to=2-3, looseness=4]
\end{tikzcd}
\end{equation}

\medskip 

\item Case $2$: 
\begin{equation}
\label{eqn:dynkin d 2}
\begin{tikzcd}
&&&& {N-1} \\
\cdots && {N-2} \\
&&&& N
\arrow[shift left=2pt, from=2-1, to=2-3]
\arrow[shift left=2pt, from=2-3, to=2-1]
\arrow["q^{\pm 1}", shift left=2pt, from=2-3, to=1-5]
\arrow["q^{\pm 1}", shift left=2pt, from=1-5, to=2-3]
\arrow["q^{\pm 1}", shift left=2pt, from=2-3, to=3-5]
\arrow["q^{\pm 1}", shift left=2pt, from=3-5, to=2-3]
\arrow["{q^{\mp 2}}", shift left=2pt, from=1-5, to=3-5]
\arrow["{q^{\mp 2}}", shift left=2pt, from=3-5, to=1-5]
\arrow[dashed, from=2-3, to=2-3, looseness=4]
\end{tikzcd}
\end{equation}

\end{itemize}

\medskip 

\noindent Some of the decorated Dynkin diagrams above also appeared in \cite{Bao}, in slightly different notation and conventions.

\medskip 

\subsection{The shuffle realization}
\label{sub:shuffle realization}

To the decorated Dynkin diagrams \eqref{eqn:dynkin b 1}-\eqref{eqn:dynkin d 2}, one can associate the zeta functions
\begin{equation} 
\label{eqn:zeta b,c,d}
\Big\{ \zeta^B_{ij}(x), \zeta^C_{ij}(x), \zeta^D_{ij}(x) \Big\}_{i,j \in I}
\end{equation} 
where $I = \{1,\dots,N\}$, according to the rule \eqref{eqn:zeta}. Using the zeta functions \eqref{eqn:zeta b,c,d}, one can define the (half) pre-quantum affine algebras
\begin{equation}
\label{eqn:quadratic b}
\tub = \BK \langle e_{i,d} \rangle_{i \in I, d\in \BZ} \Big/ \left(e_i(x)e_j(y) \zeta^{B}_{ji}\left(\frac yx\right) = e_j(y) e_i(x) \zeta^{B}_{ij}\left(\frac xy\right) \right)_{i,j \in I}
\end{equation}
\begin{equation}
\label{eqn:quadratic c}
\tuc = \BK \langle e_{i,d} \rangle_{i \in I, d\in \BZ} \Big/ \left(e_i(x)e_j(y) \zeta^{C}_{ji}\left(\frac yx\right) = e_j(y) e_i(x) \zeta^{C}_{ij}\left(\frac xy\right) \right)_{i,j \in I}
\end{equation}
\begin{equation}
\label{eqn:quadratic d}
\tud = \BK \langle e_{i,d} \rangle_{i \in I, d\in \BZ} \Big/ \left(e_i(x)e_j(y) \zeta^{D}_{ji}\left(\frac yx\right) = e_j(y) e_i(x) \zeta^{D}_{ij}\left(\frac xy\right) \right)_{i,j \in I}
\end{equation}
and the corresponding shuffle algebras $\oCS^+_B \subset \CV^+_B$, $\oCS^+_C \subset \CV^+_C$, $\oCS^+_D \subset \CV^+_D$. 

\medskip

\begin{theorem}
\label{thm:b,c,d 1}

(\cite[Theorems 1.1-1.3]{NT}) The kernels of the various homomorphisms $\widetilde{\Upsilon}^+$ are generated by the  Drinfeld-Serre relations defined in \cite{BFK,XZ}, and thus
\begin{align*}
&\ub = \tub \Big/ \emph{Ker }\widetilde{\Upsilon}^+ \\
&\uc = \tuc \Big/ \emph{Ker }\widetilde{\Upsilon}^+ \\ 
&\ud = \tud \Big/ \emph{Ker }\widetilde{\Upsilon}^+
\end{align*}
are the (Drinfeld new) halves of the quantum affine orthosymplectic superalgebras \footnote{Recall from the introduction that the phrase ``quantum affine superalgebra" of $\fosp$ type refers to the algebras defined in \cite{BFK, XZ}, which are only conjecturally isomorphic to the quantum affine superalgebras of \cite{Y2}.}.

\end{theorem}

\medskip 

\begin{theorem}
\label{thm:b,c,d 2}

(\cite[Theorems 1.1-1.3]{NT}) We have the shuffle realization
\[
\oCS^+_B = \CS^+_B, \quad \oCS^+_C = \CS^+_C, \quad \oCS^+_D = \CS^+_D
\]
where $\CS^+_B, \CS^+_C, \CS^+_D$ consist of all Laurent polynomials which vanish at the wheels \eqref{eqn:wheel finite even} and \eqref{eqn:wheel finite odd} corresponding to any two and three vertices (respectively) that span a type $A$ decorated Dynkin subdiagram \eqref{eqn:dynkin a}, together with the following additional wheels:

\medskip 

\noindent \underline{Type $B$:} 

\medskip 

\begin{itemize}[leftmargin=*]

\item \underline{Case 1:} 

\begin{equation}
\label{eqn:wheel b 1}
\begin{tikzcd}
&& N-1  && \\ \\ 
N && N && N
\arrow["{q^{\pm 1}}"', swap, from=3-1, to=1-3]
\arrow["{q^{\mp 1}}"', swap, from=3-3, to=3-1]
\arrow["{q^{\mp 1}}"', swap, from=3-5, to=3-3]
\arrow["{q^{\pm 1}}"', swap, from=1-3, to=3-5]
\end{tikzcd}
\end{equation}

\medskip 

\item \underline{Case 2:} 

\begin{equation}
\label{eqn:wheel b 2}
\begin{tikzcd}
&& N \\ \\ \\
N &&&& N
\arrow["q^{\pm 1}"', swap, from=4-1, to=1-3]
\arrow["{q^{\mp 2}}"', swap, from=4-5, to=4-1]
\arrow["q^{\pm 1}"', swap, from=1-3, to=4-5]
\end{tikzcd}
\end{equation}

\end{itemize} 

\medskip 

\noindent \underline{Type $C$:} 

\medskip 

\begin{itemize}[leftmargin=*]

\item \underline{Case 1:} 

\begin{equation}
\label{eqn:wheel c 1}
\begin{tikzcd}
&&&& N-1 \\ \\
&& N-1 \\ &&&  \\
N-1 &&&& N
\arrow["{q^{\pm 2}}"', from=5-1, to=5-5]
\arrow["{q^{\pm 2}}"', from=5-5, to=1-5]
\arrow["{q^{\mp 2}}"', from=1-5, to=3-3]
\arrow["{q^{\mp 2}}"', from=3-3, to=5-1]
\end{tikzcd}
\end{equation}

\medskip 

\item \underline{Case 2:} 

\begin{equation}
\label{eqn:wheel c 2}
\begin{tikzcd}
&& N-2 && N-1 \\ \\
N-2 && N-1 \\ &&&  \\
N-1 &&&& N
\arrow["{q^{\pm 2}}"', from=5-1, to=5-5]
\arrow["{q^{\pm 2}}"', from=5-5, to=1-5]
\arrow["{q^{\mp 1}}"', from=1-5, to=1-3]
\arrow["{q^{\mp 1}}"', from=1-3, to=3-3]
\arrow["{q^{\mp 1}}"', from=3-3, to=3-1]
\arrow["{q^{\mp 1}}"', from=3-1, to=5-1]
\end{tikzcd}
\end{equation}

\medskip 

\item \underline{Case 3:} 

\begin{equation}
\label{eqn:wheel c 3}
\begin{tikzcd}
N-3 && N-2 && N-1 \\ \\
N-2 && N-1 \\ &&& \boxed{2} \\
N-1 &&&& N
\arrow["{q^{\pm 2}}"', from=5-1, to=5-5]
\arrow["{q^{\pm 2}}"', from=5-5, to=1-5]
\arrow["{q^{\mp 1}}"', from=1-5, to=1-3]
\arrow["{q^{\mp 1}}"', swap, from=1-3, to=3-3]
\arrow["{q^{\mp 1}}"', from=3-1, to=5-1]
\arrow["{q^{\mp 1}}"', from=1-1, to=3-1]
\arrow["{q^{\mp 1}}"', from=1-3, to=1-1]
\arrow["{q^{\mp 1}}"', swap, from=3-3, to=3-1]
\arrow["{q^{\pm 2}}", from=3-1, to=1-3]
\end{tikzcd}
\end{equation}
(the meaning of the boxed $\boxed{2}$ will be explained in the next Subsection).

\medskip 

\noindent \underline{Type $D$:} 

\medskip 

\item \underline{Case 1:} no extra wheels 

\medskip 

\item \underline{Case 2:} 

\begin{equation}
\label{eqn:wheel d 2}
\begin{tikzcd}
&& N-1 &&&& N-1 \\
N-2 &&&& N-2 \\
&& N &&&& N
\arrow["{q^{\mp 2}}", from=1-3, to=3-3]
\arrow["q^{\pm 1}"', from=1-7, to=2-5]
\arrow["q^{\pm 1}", from=2-1, to=1-3]
\arrow["q^{\pm 1}"', from=2-5, to=3-7]
\arrow["q^{\pm 1}", from=3-3, to=2-1]
\arrow["{q^{\mp 2}}"', from=3-7, to=1-7]
\end{tikzcd}
\end{equation}

\medskip 

\end{itemize}

\end{theorem}

\medskip 

\subsection{Scary wheels} 
\label{sub:scary}

Let us now explain the meaning of the wheel \eqref{eqn:wheel c 3}, which is the only one in the list above which does not fall under the general framework of \eqref{eqn:wheel}. By definition, the wheel in question indicates that Laurent polynomials $E(z_{iu})_{i \in I, u \geq 1} \in \CV^+_C$ vanish when one specializes
\begin{align*}
&z_{N-3,1} \mapsto xq^{\pm 2} \\
&(z_{N-2,1}, z_{N-2,2}) \mapsto (xq^{\pm 3}, xq^{\pm 1}) \\
&(z_{N-1,1}, z_{N-1,2}, z_{N-1,3}) \mapsto (xq^{\pm 4}, xq^{\pm 2},x) \\
&z_{N,1} \mapsto xq^{\pm 2} 
\end{align*} 
However, due to the two small triangles in the corner of diagram \eqref{eqn:wheel c 3}, Laurent polynomials $E$ that satisfy the wheel conditions already vanish at the specialization above. Thus, the meaning of the $\boxed{2}$ in diagram \eqref{eqn:wheel c 3} is that we require our Laurent polynomials $E$ to moreover satisfy the following order $\geq 2$ vanishing condition
\begin{equation}
\label{eqn:rigorous}
\begin{split}
& (z_{N-2,2}-xq^{\pm 1})^2 \quad \text{divides} \\ 
& \qquad E \Big|_{(z_{N-3,1},z_{N-2,1},z_{N-1,1},z_{N-1,2},z_{N-1,3},z_{N,1}) \mapsto (xq^{\pm 2},xq^{\pm 3},xq^{\pm 4},xq^{\pm 2},x,xq^{\pm 2})}
\end{split}
\end{equation} 
Because of this more complicated algebraic definition, we call \eqref{eqn:wheel c 3} a \emph{scary wheel}, though we point out that such higher order vanishing conditions are expected to appear for more complicated decorated Dynkin diagrams (see for instance \cite{N Symm}).

\medskip 

\subsection{Quantum toroidal superalgebras in types $B,C,D$}
\label{sub:toroidal b,c,d}

Recall that the decorated Dynkin diagrams of type $B,C,D$ consist of gluing together the following:

\medskip 

\begin{itemize}[leftmargin=0.7cm]
	
\item a \emph{body}, i.e.\ a type $A$ decorated Dynkin diagram with vertices $1,\dots,N-1$ (in types $B,C$) or $1,\dots,N-2$ (in type $D$). We call the largest vertex the~\emph{neck};
	
\medskip 
	
\item a \emph{head}, i.e. one of the diagrams below
\begin{equation}
\label{eqn:type b head}
\text{Type $B$:} \qquad 
\begin{tikzcd}
N-1 && N
\arrow["q^{\pm 1}", shift left=2pt, from=1-3, to=1-1]
\arrow["q^{\pm 1}", shift left=2pt, from=1-1, to=1-3]
\arrow["q^{\mp 1}"',  from=1-3, to=1-3, looseness=4]
\arrow[dashed, from=1-1, to=1-1, looseness=4]
\end{tikzcd} \quad \begin{tikzcd}
N-1 && N 
\arrow["q^{\pm 1}", shift left=2pt, from=1-3, to=1-1]
\arrow["q^{\pm 1}", shift left=2pt, from=1-1, to=1-3]
\arrow["q^{\mp 2}"', from=1-3, to=1-3, out=45, in=135, distance=10mm]
\arrow["q^{\pm 1}", from=1-3, to=1-3, out=315, in=225, distance=10mm]
\arrow[dashed, from=1-1, to=1-1, looseness=4]
\end{tikzcd}
\end{equation}
\begin{equation}
\label{eqn:type c head}
\text{Type $C$:} \qquad \qquad
\begin{tikzcd}
N-1 && N
\arrow["q^{\pm 2}", shift left=2pt, from=1-3, to=1-1]
\arrow["q^{\pm 2}", shift left=2pt, from=1-1, to=1-3]
\arrow["q^{\mp 4}"', from=1-3, to=1-3, looseness=4]
\arrow[dashed, from=1-1, to=1-1, looseness=4]
\end{tikzcd} 
\end{equation}
\begin{equation}
\label{eqn:type d head}
\text{Type $D$:} \quad 
\begin{tikzcd}
&& {N-1} \\
{N-2} \\
&& N
\arrow["q^{\pm 1}", shift left=2pt, from=2-1, to=1-3]
\arrow["q^{\pm 1}", shift left=2pt, from=1-3, to=2-1]
\arrow["q^{\pm 1}", shift left=2pt, from=2-1, to=3-3]
\arrow["q^{\pm 1}", shift left=2pt, from=3-3, to=2-1]
\arrow["q^{\mp 2}", swap, from=1-3, to=1-3, looseness=4]
\arrow["q^{\mp 2}", swap, from=3-3, to=3-3, looseness=4]
\arrow[dashed, from=2-1, to=2-1, looseness=4]
\end{tikzcd} \quad \begin{tikzcd}
&& {N-1} \\
{N-2} \\
&& N
\arrow["q^{\pm 1}", shift left=2pt, from=2-1, to=1-3]
\arrow["q^{\pm 1}", shift left=2pt, from=1-3, to=2-1]
\arrow["q^{\pm 1}", shift left=2pt, from=2-1, to=3-3]
\arrow["q^{\pm 1}", shift left=2pt, from=3-3, to=2-1]
\arrow["{q^{\mp 2}}", shift left=2pt, from=1-3, to=3-3]
\arrow["{q^{\mp 2}}", shift left=2pt, from=3-3, to=1-3]
\arrow[dashed, from=2-1, to=2-1, looseness=4]
\end{tikzcd}
\end{equation}

\end{itemize}

\medskip 
\noindent (all the signs $\pm$ in a given diagram are one and the same). Similarly, an \emph{extended decorated Dynkin diagram} of type $B,C,D$ consists of a type $A$ decorated Dynkin diagram (a body) to which we glue two heads at the leftmost and rightmost ends, respectively. An extended decorated Dynkin diagram has two \emph{neck} vertices where the body is connected to the two heads, and thus two neck parameters $q^{\pm 1}$ and $q^{\pm' 1}$ connecting the respective necks with the outermost vertices. The vertex set of an extended decorated Dynkin diagram will be denoted by $\tilde{I} = \{0,\dots,N\}$, and the following are the possible extended decorated Dynkin diagrams:
\begin{equation}
\label{eqn:dynkin bb}
\text{Type $BB$:} \begin{tikzcd}
0 && 1 & \cdots & {N-1} && N
\arrow["q^{\pm 1}", shift left=2pt, from=1-7, to=1-5]
\arrow["q^{\pm 1}", shift left=2pt, from=1-5, to=1-7]
\arrow[shift left=2pt, from=1-5, to=1-4]
\arrow[shift left=2pt, from=1-4, to=1-5]
\arrow[shift left=2pt, from=1-4, to=1-3]
\arrow[shift left=2pt, from=1-3, to=1-4]
\arrow["q^{\pm' 1}", shift left=2pt, from=1-3, to=1-1]
\arrow["q^{\pm' 1}", shift left=2pt, from=1-1, to=1-3]
\arrow["q^{\mp 1}"', from=1-7, to=1-7, looseness=3.3]
\arrow[dashed, from=1-5, to=1-5, looseness=3.3]
\arrow[dashed, from=1-3, to=1-3, looseness=4.5]
\arrow["q^{\mp' 1}"', from=1-1, to=1-1, looseness=4.5]
\end{tikzcd}
\end{equation}
\begin{equation}
\label{eqn:dynkin bc}
\text{Type $BC$:} \begin{tikzcd}
0 && 1 & \cdots & {N-1} && N
\arrow["q^{\pm 2}", shift left=2pt, from=1-7, to=1-5]
\arrow["q^{\pm 2}", shift left=2pt, from=1-5, to=1-7]
\arrow[shift left=2pt, from=1-5, to=1-4]
\arrow[shift left=2pt, from=1-4, to=1-5]
\arrow[shift left=2pt, from=1-4, to=1-3]
\arrow[shift left=2pt, from=1-3, to=1-4]
\arrow["q^{\pm' 1}", shift left=2pt, from=1-3, to=1-1]
\arrow["q^{\pm' 1}", shift left=2pt, from=1-1, to=1-3]
\arrow["q^{\mp 4}"', from=1-7, to=1-7, looseness=3.3]
\arrow[dashed, from=1-5, to=1-5, looseness=3.3]
\arrow[dashed, from=1-3, to=1-3, looseness=4.5]
\arrow["q^{\mp' 1}"', from=1-1, to=1-1, looseness=4.5]
\end{tikzcd}
\end{equation}
\begin{equation}
\label{eqn:dynkin bd}
\text{Type $BD$:} \quad \begin{tikzcd}
&&&&& N-1 \\
0 && 1 & \cdots & {N-2} \\
&&&&& N
\arrow[shift left=2pt, from=2-5, to=2-4]
\arrow[shift left=2pt, from=2-4, to=2-5]
\arrow[shift left=2pt, from=2-4, to=2-3]
\arrow[shift left=2pt, from=2-3, to=2-4]
\arrow["q^{\pm' 1}", shift left=2pt, from=2-3, to=2-1]
\arrow["q^{\pm' 1}", shift left=2pt, from=2-1, to=2-3]
\arrow["q^{\pm 1}", shift left=2pt, from=1-6, to=2-5]
\arrow["q^{\pm 1}", shift left=2pt, from=2-5, to=1-6]
\arrow["q^{\pm 1}", shift left=2pt, from=3-6, to=2-5]
\arrow["q^{\pm 1}", shift left=2pt, from=2-5, to=3-6]
\arrow["q^{\mp 2}", shift left=2pt, from=3-6, to=1-6]
\arrow["q^{\mp 2}", shift left=2pt, from=1-6, to=3-6]
\arrow["q^{\mp' 1}"', from=2-1, to=2-1, looseness=4.5]
\arrow[dashed, from=2-3, to=2-3, looseness=4.5]
\arrow[dashed, from=2-5, to=2-5, looseness=3.3]
\end{tikzcd}
\end{equation}
\begin{equation}
\label{eqn:dynkin cc}
\text{Type $CC$:} \begin{tikzcd}
0 && 1 & \cdots & {N-1} && N
\arrow["q^{\pm 2}", shift left=2pt, from=1-7, to=1-5]
\arrow["q^{\pm 2}", shift left=2pt, from=1-5, to=1-7]
\arrow[shift left=2pt, from=1-5, to=1-4]
\arrow[shift left=2pt, from=1-4, to=1-5]
\arrow[shift left=2pt, from=1-4, to=1-3]
\arrow[shift left=2pt, from=1-3, to=1-4]
\arrow["q^{\pm' 2}", shift left=2pt, from=1-3, to=1-1]
\arrow["q^{\pm' 2}", shift left=2pt, from=1-1, to=1-3]
\arrow["q^{\mp 4}"', from=1-7, to=1-7, looseness=3.3]
\arrow[dashed, from=1-5, to=1-5, looseness=3.3]
\arrow[dashed, from=1-3, to=1-3, looseness=4]
\arrow["q^{\mp' 4}"', from=1-1, to=1-1, looseness=4]
\end{tikzcd}
\end{equation}
\begin{equation}
\label{eqn:dynkin cd}
\text{Type $CD$:} \quad \begin{tikzcd}
&&&&& N-1 \\
0 && 1 & \cdots & {N-2} \\
&&&&& N
\arrow[shift left=2pt, from=2-5, to=2-4]
\arrow[shift left=2pt, from=2-4, to=2-5]
\arrow[shift left=2pt, from=2-4, to=2-3]
\arrow[shift left=2pt, from=2-3, to=2-4]
\arrow["q^{\pm' 2}", shift left=2pt, from=2-3, to=2-1]
\arrow["q^{\pm' 2}", shift left=2pt, from=2-1, to=2-3]
\arrow["q^{\pm 1}", shift left=2pt, from=1-6, to=2-5]
\arrow["q^{\pm 1}", shift left=2pt, from=2-5, to=1-6]
\arrow["q^{\pm 1}", shift left=2pt, from=3-6, to=2-5]
\arrow["q^{\pm 1}", shift left=2pt, from=2-5, to=3-6]
\arrow["q^{\mp 2}", shift left=2pt, from=3-6, to=1-6]
\arrow["q^{\mp 2}", shift left=2pt, from=1-6, to=3-6]
\arrow["q^{\mp' 4}"', from=2-1, to=2-1, looseness=4.5]
\arrow[dashed, from=2-3, to=2-3, looseness=4.5]
\arrow[dashed, from=2-5, to=2-5, looseness=3.3]
\end{tikzcd}
\end{equation}
\begin{equation}
\label{eqn:dynkin dd}
\text{Type $DD$:} \quad \begin{tikzcd}
0 &[1cm]&&& N-1 \\
& 2 & \cdots & {N-2} & \\
1 &&&& N
\arrow["q^{\pm' 1}", shift left=2pt, from=1-1, to=2-2]
\arrow["q^{\pm' 1}", shift left=2pt, from=2-2, to=1-1]
\arrow["q^{\pm' 1}", shift left=2pt, from=3-1, to=2-2]
\arrow["q^{\pm' 1}", shift left=2pt, from=2-2, to=3-1]
\arrow["q^{\pm 1}", shift left=2pt, from=1-5, to=2-4]
\arrow["q^{\pm 1}", shift left=2pt, from=2-4, to=1-5]
\arrow["q^{\pm 1}", shift left=2pt, from=3-5, to=2-4]
\arrow["q^{\pm 1}", shift left=2pt, from=2-4, to=3-5]
\arrow["q^{\mp 2}", shift left=2pt, from=3-5, to=1-5]
\arrow["q^{\mp 2}", shift left=2pt, from=1-5, to=3-5]
\arrow["q^{\mp' 2}", shift left=2pt, from=3-1, to=1-1]
\arrow["q^{\mp' 2}", shift left=2pt, from=1-1, to=3-1]
\arrow[shift left=2pt, from=2-4, to=2-3]
\arrow[shift left=2pt, from=2-3, to=2-4]
\arrow[shift left=2pt, from=2-3, to=2-2]
\arrow[shift left=2pt, from=2-2, to=2-3]
\arrow[dashed,from=2-2, to=2-2, looseness=4.5]
\arrow[dashed,  from=2-4, to=2-4, looseness=3.3]
\end{tikzcd}
\end{equation}
(in types $B$ and $D$ we may replace either of the heads represented above by the other head in \eqref{eqn:type b head}, \eqref{eqn:type d head}). The signs $\pm'$ and $\pm$ are independent of each other. One also has types $CB$, $DB$ and $DC$ extended Dynkin diagrams by reflecting the above diagrams left-to-right. Any choice of a root system of a type $B,C,D$ affine Lie superalgebra yields a extended decorated Dynkin diagram (see \cite{Y2}), although the converse is not true. We henceforth assume that the two necks are different (and moreover non-adjacent in types $CC$), which eliminates some small rank cases.

\medskip 

\begin{definition}
\label{def:special}

Let $Q \in q^{\BZ}$ be the product of all the horizontal edge parameters in an extended decorated Dynkin diagram (for each type $D$ head, we include the corresponding neck parameter exactly once in the product), times a single factor of $q^{\pm 1}$ in those types $BB$ where both heads have a loop parameter of $q^{\pm 1}$. An extended decorated Dynkin diagram is called \emph{special} if it is of types $BB, CC, DD$, $CD, DC$ and $Q = 1$, and \emph{non-special} otherwise.

\end{definition}

\medskip

\noindent Let $\widetilde{U}^+_q(L\fg_\Gamma)$ denote the algebra \eqref{eqn:pre} corresponding the the zeta functions encoded by an extended decorated Dynkin diagram $\Gamma$ as above.

\medskip

\begin{theorem}
\label{thm:toroidal b,c,d 1}
	
(\cite[Theorem 1.5]{NT}) For any non-special extended Dynkin diagram $\Gamma$, $\emph{Ker }\widetilde{\Upsilon}^+$ is generated by the relations defined in \cite{BFK,XZ} corresponding to the two heads and the relations defined in \cite{Y2} corresponding to the body. Therefore, we call
\[
U^+_q(L\fg_\Gamma) = \widetilde{U}^+_q(L\fg_\Gamma) \Big/ \emph{Ker }\widetilde{\Upsilon}^+
\]
a quantum toroidal Lie superalgebra of type $XY$, where $X,Y \in \{B,C,D\}$ denote the shape of the two heads of $\Gamma$ (though strictly speaking, the name ``toroidal" should only be applied to those $\Gamma$ that come from affinizing finite type Lie superalgebras).
	
\end{theorem}

\medskip

\begin{theorem}
\label{thm:toroidal b,c,d 2}

(\cite{NT}) For any non-special extended Dynkin diagram $\Gamma$ (aside from the small rank cases referenced in the sentence preceding Definition~\ref{def:special}) we have
\[
\oCS^+_{\Gamma} = \CS_{\Gamma}^+
\]
where $\CS^+_{\Gamma}$ consists of all Laurent polynomials which vanish at the wheels \eqref{eqn:wheel finite even}-\eqref{eqn:wheel finite odd} corresponding to the body of $\Gamma$ and the wheels \eqref{eqn:wheel b 1}-\eqref{eqn:wheel d 2} corresponding to the two heads of $\Gamma$.

\end{theorem}

\medskip 

\subsection{Quantum toroidal superalgebras in type $A$}
\label{sub:toroidal a}

Theorems~\ref{thm:toroidal b,c,d 1} and~\ref{thm:toroidal b,c,d 2} have analogues in type $A$. Let us take
\[
\fg_{\widehat{A}} = \widehat{\fsl}_{m|n}
\]
defined with respect to any parity sequence \eqref{eqn:sequence}. The corresponding decorated Dynkin diagram is the following, with vertex set 
\[
\widehat{I} = \{0,\dots,N-1\} = \BZ/N\BZ
\] 
and parameters in the field $\BQ(q,d)$:

\begin{equation}
\label{eqn:dynkin affine a}
\begin{tikzcd}[row sep=3em,column sep=3em]
&& {N-1} &&& {0} && \\ \\
\vdots &&&&&&& {1} \\ \\
&& {3} &&& {2}
\arrow["{q^{-s_N}d^{-s_N}}", shift left=2pt, from=1-3, to=1-6]
\arrow["{q^{-s_N}d^{s_N}}", shift left=2pt, from=1-6, to=1-3]
\arrow["{q^{-s_1}d^{-s_1}}", shift left=2pt, from=1-6, to=3-8]
\arrow["{q^{-s_1}d^{s_1}}", shift left=2pt, from=3-8, to=1-6]
\arrow[shift left=2pt, from=3-1, to=1-3]
\arrow[shift left=2pt, from=1-3, to=3-1]
\arrow["{q^{-s_2}d^{-s_2}}", shift left=2pt, from=3-8, to=5-6]
\arrow["{q^{-s_2}d^{s_2}}", shift left=2pt, from=5-6, to=3-8]
\arrow[shift left=2pt, from=5-3, to=3-1]
\arrow[shift left=2pt, from=3-1, to=5-3]
\arrow["{q^{-s_3}d^{-s_3}}", shift left=2pt, from=5-6, to=5-3]
\arrow["{q^{-s_3}d^{s_3}}",shift left=2pt, from=5-3, to=5-6]
\arrow["q^{2s_{N-1}}"', distance=4em, dashed, from=1-3, to=1-3, out=45, in=135]
\arrow["q^{2s_{N}}"', distance=4em, dashed, from=1-6, to=1-6, out=45, in=135]
\arrow["q^{2s_1}", distance=4.5em, dashed, from=3-8, to=3-8, out=45, in=315]
\arrow["q^{2s_2}", distance=4.5em, dashed, from=5-6, to=5-6, out=315, in=225]
\arrow["q^{2s_{3}}", distance=4.5em, dashed, from=5-3, to=5-3, out=315, in=225]
\end{tikzcd}
\end{equation}
We henceforth assume that $m\neq n$, which is analogous to the non-special condition of Definition~\ref{def:special}. If one applies the discussion of Subsection~\ref{sub:algebras} to the zeta functions encoded by the quiver \eqref{eqn:dynkin affine a}, one obtains quantum toroidal superalgebras
\[
\tuwa \twoheadrightarrow \uwa 
\]
and shuffle algebras 
\[
\oCS^+_{\widehat{A}} \subset \CV^+_{\widehat{A}}
\]
We have the following analogues of Theorems~\ref{thm:a 1} and~\ref{thm:a 2}.

\medskip 

\begin{theorem}
\label{thm:wa 1}
	
(\cite{NT}) Assuming $N \geq 4$ and $m \neq n$, the kernel of the map $\widetilde{\Upsilon}^+$ is generated by the cubic and quartic Drinfeld-Serre relations defined in \cite{BM}, and so
\[
\uwa = \tuwa \Big/ \emph{Ker } \widetilde{\Upsilon}^+
\]
is the (Drinfeld new like) half of the type $A$ quantum toroidal superalgebra. 
	
\end{theorem}

\medskip 

\begin{theorem}
\label{thm:wa 2}
	
(\cite{N Red}) We have
\[
\oCS^+_{\widehat{A}} = \CS^+_{\widehat{A}}
\]
where $\CS^+_{\widehat{A}} \subset \CV^+_{\widehat{A}}$ consists of all Laurent polynomials which vanish at the wheels
\begin{equation}
\label{eqn:wheel affine even}
\begin{tikzcd}
i-1 && i &&&& i \\ \\
i &&&& i && i+1
\arrow["{q^{-s_{i}}d^{-s_{i}}}"', from=1-1, to=3-1]
\arrow["{q^{-s_i}d^{s_i}}"', from=1-3, to=1-1]
\arrow["{q^{-s_i}d^{-s_i}}", from=1-7, to=3-7]
\arrow["{q^{2s_i}}"', from=3-1, to=1-3]
\arrow["{q^{2s_i}}", from=3-5, to=1-7]
\arrow["{q^{-s_i}d^{s_i}}", from=3-7, to=3-5]
\end{tikzcd}	
\end{equation}
for any $i \in \hat{I}$ such that $s_i = s_{i+1}$, and
\begin{equation}
\label{eqn:wheel affine odd}
\begin{tikzcd}
i-1 && i \\ \\
i && i+1
\arrow["{q^{-s_i}d^{-s_i}}", from=1-1, to=1-3]
\arrow["{q^{-s_{i+1}}d^{-s_{i+1}}}", from=1-3, to=3-3]
\arrow["{q^{-s_i}d^{s_i}}", from=3-1, to=1-1]
\arrow["{q^{-s_{i+1}}d^{s_{i+1}}}", from=3-3, to=3-1]
\end{tikzcd}
\end{equation}
for any $i \in \hat{I}$ such that $s_i \neq s_{i+1}$ (we define the parity sequence $s_i$ for all $i \in \BZ$ by using a choice \eqref{eqn:sequence} and the periodicity $s_{i+N} = s_i$ for all $i \in \BZ$.)
\end{theorem}

\bigskip

\section{$K$-theoretic Hall algebras of quivers with potential}
\label{sec:geometry}

\medskip 

\subsection{General definitions}
\label{sub:general}

Let us consider a quiver $Q$ with vertex set $I$ and arrow set $\Omega$. We will also fix a potential $W$, i.e. a linear combination of oriented cycles of $Q$ modulo cyclic permutations. Consider an algebraic torus $T$ endowed with characters
\[
\Big\{ t_\alpha : T \rightarrow \BC^* \Big\}_{\alpha \in \Omega}
\]
We make the important assumption that for every summand which appears with non-zero coefficient in the potential
\begin{equation}
\label{eqn:potential}
W = \dots + (\text{non-zero constant}) \cdot \alpha_1 \dots \alpha_n + \cdots
\end{equation}
we have 
\begin{equation}
\label{eqn:assumption}
t_{\alpha_1}\dots t_{\alpha_n} = 1
\end{equation}
as well as the technical assumption that $t_\alpha$ is non-trivial for any loop $\alpha$. For any $\bn = (n_i \geq 0)_{i \in I}$, we consider the affine space of $\bn$-dimensional representations of $Q$:
\begin{equation}
\label{eqn:quiver representation}
\text{Rep}_{\bn}(Q) = \Big\{ \BC^{n_i} \xrightarrow{\phi_\alpha} \BC^{n_j} \Big\}_{\forall \alpha : i \rightarrow j \text{ in } \Omega}
\end{equation}
Let $G_{\bn} = \prod_{i \in I} GL_{n_i}(\BC)$. We have the commuting actions 
\[
G_{\bn} \curvearrowright \text{Rep}_{\bn}(Q) \curvearrowleft T
\]
where $G_{\bn}$ acts by conjugation on collections of linear maps \eqref{eqn:quiver representation}, and $T$ scales each $\phi_\alpha$ by the character $t_\alpha$. Because of \eqref{eqn:assumption}, the trace of the potential induces a $T$-invariant function on the representation stack
\begin{equation}
\label{eqn:stack and function}
\text{Rep}_{\bn}(Q)/G_{\bn} \xrightarrow{\text{tr}(W)} \BC 
\end{equation}
simply by sending \eqref{eqn:potential} to $\cdots + (\text{non-zero constant}) \cdot \text{tr} \left(\phi_{\alpha_1} \cdots \phi_{\alpha_n}\right) + \cdots$.

\medskip 

\subsection{Matrix factorizations and critical $K$-theory}
\label{sub:mf}

Given a smooth quotient stack $X$ over $\BC$, acted on by a torus $T$, and endowed with a $T$-invariant function 
\[
f : X \rightarrow \BC
\]
one can define the dg enhancement of the absolute derived category of $T$-equivariant coherent factorizations (see \cite{P2} for foundations on the contents of this Subsection)
\[
\text{MF}^T(X,f)
\]
The Grothendieck group of this category is sometimes called the \emph{critical $K$-theory}
\[
K_{\text{crit}}^T(X,f) = K_0 ( \text{MF}^T(X,f) )
\]
When $X = \text{Rep}_{\bn}(Q)/G_{\bn}$ and $f = \text{tr}(W)$ are as in \eqref{eqn:stack and function}, the direct sum 
\begin{equation}
\label{eqn:k-ha}
K_{\text{crit}}^T(Q,W) = \bigoplus_{\bn = (n_i \geq 0)_{i \in I}} K_{\text{crit}}^T \left( \text{Rep}_{\bn}(Q)/G_{\bn}, \text{tr}(W) \right) 
\end{equation}
was made into an algebra called the $K$-theoretic Hall algebra (\cite{P}, inspired by \cite{KS}) using a certain convolution diagram. We note that one has the freedom to insert various line bundle twists in this convolution diagram, the choices of which correspond to various monomials concealed by the $\sim$ symbol in \eqref{eqn:zeta}. Let
\[
\BK = \text{Frac}(\text{Rep}_T) = \text{Frac} \left(K^T(\text{point})\right) 
\]
and we will consider the localized $K$-HA
\[
K_{\text{crit}}^T(Q,W)_{\text{loc}} = K_{\text{crit}}^T(Q,W) \bigotimes_{\text{Rep}_T} \BK 
\]
We will not need to review the algebra structure on \eqref{eqn:k-ha}, and list here the only two properties we will use. Firstly, when $W = 0$, we have natural identifications
\[
K_{\text{crit}}^T(Q,0)_{\text{loc}} = \bigoplus_{\bn} K^{G_{\bn} \times T} \left( \text{Rep}_{\bn}(Q) \right)_{\text{loc}} = \bigoplus_{\bn} K^{G_{\bn} \times T} \left( \text{point} \right)_{\text{loc}} = \CV^+
\]
with the multiplication on the LHS corresponding to the shuffle product \eqref{eqn:shuffle product} on the RHS. Secondly, we have an algebra homomorphism
\[
K_{\text{crit}}^T(Q,W)_{\text{loc}} \longrightarrow K_{\text{crit}}^T(Q,0)_{\text{loc}}
\]
whose image is contained in the image of $\jmath_*$ (\cite[Proposition 3.6]{P2}, whose hypothesis precisely requires our assumption that $t_\alpha$ is non-trivial for any loop $\alpha$), where 
\[
\jmath : \{\text{tr}(W) = 0\} \hookrightarrow \bigsqcup_{\bn = (n_i\geq 0)_{i \in I}} \text{Rep}_{\bn}(Q)/G_{\bn}
\]
is the closed embedding of the zero locus of the function $\text{tr}(W)$. Putting everything together, we have a commutative diagram 
\begin{equation}
\label{eqn:commutative diagram}
\begin{tikzcd}
K_{\text{crit}}^T(Q,W)_{\text{loc}} \arrow[rr, "\iota"] \arrow[dr, dashed] && \CV^+ \\
& \text{Im } \jmath_* \arrow[ur, hook] &
\end{tikzcd}
\end{equation}

\medskip 

\subsection{Quivers with potential for Lie superalgebras}
\label{sub:quiver with potential}

Quivers with potential associated to quantum affine (super)algebras have been intensively studied in the mathematical physics literature (\cite{GLY, GY, LY, NW}, and particularly \cite{Bao} for many of the orthosymplectic types considered in the present paper). However, the discussion below will indicate a perhaps unexpected quiver with potential that we associate to the type $C$ Case 3 Dynkin diagram \eqref{eqn:dynkin c 3}, and the proof of Theorem~\ref{thm:main} will explain why this particular choice of quiver is necessary. For any Lie superalgebra
\[
\fg \in \Big\{\fsl_{m|n}, \fosp_{m|n}, \widehat{\fosp}_{m|n}, \widehat{\fsl}_{m|n} \Big\}
\]
let $Q$ denote the decorated Dynkin diagram \eqref{eqn:dynkin a}, \eqref{eqn:dynkin b 1}-\eqref{eqn:dynkin d 2}, \eqref{eqn:dynkin bb}-\eqref{eqn:dynkin dd} (assume non-special), \eqref{eqn:dynkin affine a} (assume $m \neq n$), respectively. In fact, we could assume $Q$ to be an extended Dynkin diagram even if it doesn't come from an affine orthosymplectic Lie superalgebra. In each of these cases, we consider the torus
\[
T = \BC^* \quad \text{with} \quad \text{Rep}_T = \BZ[q^{\pm 1}] \quad \Rightarrow \quad \BK = \BQ(q)
\]
acting on the arrows with the weights indicated in diagrams \eqref{eqn:dynkin a}, \eqref{eqn:dynkin b 1}-\eqref{eqn:dynkin d 2}, \eqref{eqn:dynkin bb}-\eqref{eqn:dynkin dd}, respectively. In affine type $A$, we may consider instead
\[
T = \BC^* \times \BC^* \quad \text{with} \quad \text{Rep}_T = \BZ[q^{\pm 1},d^{\pm 1}] \quad \Rightarrow \quad \BK = \BQ(q,d)
\]
acting on the arrows with the weights from diagram \eqref{eqn:dynkin affine a}. In all types, we define 
\begin{equation}
\label{eqn:super potential}
W = \sum_{\text{wheels }\bigcirc \text{ as in \eqref{eqn:wheel}}} \gamma_{\bigcirc} \cdot \alpha_k \dots \alpha_1
\end{equation}
for arbitrary non-zero constants $\gamma_{\bigcirc} \in \BC^*$, where the sum goes over all applicable wheels among \eqref{eqn:wheel finite even}-\eqref{eqn:wheel finite odd}, \eqref{eqn:wheel b 1}-\eqref{eqn:wheel d 2}, \eqref{eqn:wheel affine even}-\eqref{eqn:wheel affine odd}; since the terms of a potential are only defined up to cyclic permutations, we require that the total coefficient of any given cycle (irrespective of how it is ordered as $\alpha_k \dots \alpha_1$) be non-zero.

\medskip 

\noindent The only exception to formula \eqref{eqn:super potential} occurs in ``type $C$ Case $3$" (and extended Dynkin diagrams with heads of this type) in which case we posit that the potential should include the following terms that correspond to the scary wheel \eqref{eqn:wheel c 3}:
\begin{equation}
\label{eqn:scary potential}
W = \cdots + \gamma_1 \left( \begin{tikzcd}[row sep=2em, column sep=2em]
N-3 & N-2 \\ 
N-2 
\arrow[from=1-1, to=2-1]
\arrow[from=1-2, to=1-1]
\arrow[from=2-1, to=1-2]
\end{tikzcd} \right) + \gamma_2 \left(\begin{tikzcd}[row sep=2em, column sep=2em]
 & N-2 \\ 
N-2 & N-1
\arrow[from=1-2, to=2-2]
\arrow[from=2-2, to=2-1]
\arrow[from=2-1, to=1-2]
\end{tikzcd} \right) +
\end{equation}
\[
+ \gamma_3 \left( \begin{tikzcd}[row sep=1.5em, column sep=1.5em]
N-3 & N-2 & N-1 \\ 
N-2 & \\
N-1 && N
\arrow[from=3-1, to=3-3]
\arrow[from=3-3, to=1-3]
\arrow[from=1-3, to=1-2]
\arrow[from=2-1, to=3-1]
\arrow[from=1-1, to=2-1]
\arrow[from=1-2, to=1-1]
\end{tikzcd} \right)
+ \gamma_4 \left(\begin{tikzcd}[row sep=1.5em, column sep=1.5em]
& N-2 & N-1 \\ 
N-2 & N-1 \\
N-1 && N
\arrow[from=3-1, to=3-3]
\arrow[from=3-3, to=1-3]
\arrow[from=1-3, to=1-2]
\arrow[from=1-2, to=2-2]
\arrow[from=2-1, to=3-1]
\arrow[from=2-2, to=2-1]
\end{tikzcd} \right) + \cdots
\]
where $\gamma_1,\gamma_2,\gamma_3,\gamma_4 \in \BC$ are such that $\gamma_1 \neq 0 \neq \gamma_2$ and $\gamma_1\gamma_4 \neq \gamma_2\gamma_3$. The ellipses in \eqref{eqn:scary potential} indicate terms as in \eqref{eqn:super potential} that correspond to wheels which are not depicted in \eqref{eqn:wheel c 3}.

\medskip 

\begin{remark}

In \cite{Bao}, it was postulated that the summands of $W$ should be read off the higher order relations that hold in the quantum affine superalgebra. We slightly amend this principle by positing that the summands of $W$ should be chosen such that they enforce the wheel conditions (see the proof of Theorem~\ref{thm:main} below), while the wheel conditions themselves are in one-to-one correspondence with the higher order relations that hold in the quantum affine superalgebra (see \cite{N Arb}). For example, the fourth term which appears in \eqref{eqn:scary potential} was already pointed out in \cite{Bao} in type $C$ Case $2$, but our description above shows that the same term appears in type $C$ Case $3$.
	
\end{remark}

\medskip 

\begin{proof} \emph{of Theorem~\ref{thm:main}:} By Theorems~\ref{thm:a 1}, \ref{thm:b,c,d 1}, \ref{thm:toroidal b,c,d 1}, \ref{thm:wa 1}, it suffices to show that the image of the map $\iota$ in \eqref{eqn:commutative diagram} satisfies
\begin{equation}
\label{eqn:1}
\text{Im } \iota = \oCS^+
\end{equation}
The inclusion $\supseteq$ is an immediate consequence of the defining fact that $\oCS^+$ is generated by elements of degree $\bn = (0,\dots,0,1,0,\dots,0)$ where the $1$ is located at various positions $i \in I$, and the fact that corresponding elements also lie in $K_{\text{crit}}(Q,W)_{\text{loc}}$. For the inclusion $\subseteq$, it suffices to prove that
\begin{equation}
\label{eqn:2}
\text{Im }\iota \subseteq \CS^+ 
\end{equation}
due to Theorems~\ref{thm:a 2},~\ref{thm:b,c,d 2},~\ref{thm:toroidal b,c,d 2},~\ref{thm:wa 2}. In order to prove \eqref{eqn:2}, we adapt the approach pioneered in \cite{Z}. Since the diagram \eqref{eqn:commutative diagram} commutes, it suffices to show that
\begin{equation}
\label{eqn:3}
\text{Im }\jmath_* \subseteq \CS^+ 
\end{equation} 
In other words, we take any element $E \in \text{Im }\jmath_*$ and we must show that $E$ satisfies the wheel conditions that define $\CS^+$. Fix any wheel $\bigcirc$ of the form \eqref{eqn:wheel}, and consider the following argument. For any $\bn = (n_i \geq 0)_{i \in I}$, consider the locally closed subset
\[
L \subset \text{Rep}_{\bn}(Q)
\]
consisting of quiver representations whose only non-zero arrow maps satisfy
\[
\phi_{\alpha_{s+1}}(v_{i_s,u_s}) \in \BC^* v_{i_{s+1},u_{s+1}}
\]
for all $s \in \BZ/k\BZ$ and $u_s \in \{1,\dots,n_{i_s}\}$ such that $u_s \neq u_{s'}$ if $s \neq s'$ and $i_s = i_{s'}$ (with respect to a fixed choice of basis $v_{i,1},\dots,v_{i,n_i}$ of each $\BC^{n_i}$). We have 
\[
L \cong (\BC^*)^k \qquad \text{and} \qquad L \cap \text{Im }\jmath = \varnothing
\]
Moreover, the $K$-theoretic restriction map from $\text{Rep}_{\bn}(Q)$ to $L$ precisely corresponds to specializing Laurent polynomials as in the left-hand side of \eqref{eqn:satisfies wheel}. Combining all these remarks shows that
\begin{equation}
\label{eqn:chain}
E \in \text{Im }\jmath_* \quad \Rightarrow \quad E|_L = 0 \quad \Rightarrow \quad E \text{ satisfies } \eqref{eqn:satisfies wheel}
\end{equation}
We skip the proof of the assertions above, since it closely follows that of \cite[Proposition 2.11]{N Wheel} and \cite[Corollary 1.18]{N Red}. Instead, let us now consider the special case of the Lie superalgebra of type $C$ Case $3$, and focus on the four terms displayed in \eqref{eqn:scary potential}. Our goal is to show that any element $E \in \text{Im }\jmath_*$ satisfies the scary wheel condition \eqref{eqn:rigorous}. To this end, consider the locally closed subset 
\[
L \stackrel{\kappa}\hookrightarrow \text{Rep}_{\bn}(Q)
\]
consisting of quiver representations whose only non-zero arrow maps (for brevity we will write $\phi_{ij}$ for the linear map corresponding to the unique arrow from vertex $N-i$ to $N-j$, and $v_{i,1},v_{i,2},\dots$ for a choice of basis vectors of the vector space located at vertex $N-i$) satisfy
\begin{equation}
\label{eqn:phi}
\begin{split} 
&\phi_{10}(v_{1,3}) = x_1 v_{0,1} \\
&\phi_{01}(v_{0,1}) = x_2 v_{1,1} \\ 
&\phi_{12}(v_{1,1}) = x_3 v_{2,1} \\
&\phi_{21}(v_{2,1}) = x_4 v_{1,2} \\
&\phi_{23}(v_{2,1}) = x_5 v_{3,1} \\
&\phi_{22}(v_{2,2}) = a v_{2,1} \\
&\phi_{32}(v_{3,1}) = b v_{2,2} \\
&\phi_{12}(v_{1,2}) = c v_{2,2} \\
&\phi_{21}(v_{2,2}) = d v_{1,3} 
\end{split} 
\end{equation}
where $x_1,\dots,x_5 \in \BC^*$ and $a,b,c,d \in \BC$ are required to satisfy
\begin{equation}
\label{eqn:non-zero}
\gamma_1 a'b' + \gamma_2 a'c' + \gamma_3 b'd' + \gamma_4 c'd' \neq 0
\end{equation}
where $a' = a$, $b' = b x_5$, $c' = cx_4$, $d' = dx_1x_2x_3$. The reason for the choice \eqref{eqn:non-zero} is that the left-hand side of the equation above is precisely the trace of the potential \eqref{eqn:scary potential} on the locus $L$. As a consequence, we have
\begin{equation}
\label{eqn:choice 2}
L \cap \text{Im }\jmath = \varnothing
\end{equation}
However,
\begin{equation}
\label{eqn:choice 1}
L \cong (\BC^*)^5 \times U 
\end{equation} 
where $U \subset \BA^4$ is given by \eqref{eqn:non-zero}, i.e. the complement of a quadratic hypersurface in the affine space with coordinates $a',b',c',d'$, which is smooth away from the origin (the smoothness follows from the condition $\gamma_1\gamma_4 \neq \gamma_2\gamma_3$). Consider the maximal torus
\[
H \subset G_{\bn}
\]
which acts diagonally in the basis $v_{i,1},v_{i,2},\dots$ and let us denote the corresponding elementary characters of $H$ by $z_{i,1},z_{i,2},\dots$. We have a diagram
\[
\begin{tikzcd}
G^T \left(\{\text{tr}(W) = 0\}/H \right) \arrow[d] \arrow[r,"\bar{\jmath}^*"]  & K^T \left(\text{Rep}_{\bn}(Q)/H \right) \arrow[d] \arrow[r,"\kappa^*"] & K^T \left(L/H \right) \\ 
G^T \left(\{\text{tr}(W) = 0\}/G_{\bn} \right) \arrow[r,"\jmath^*"] & K^T \left(\text{Rep}_{\bn}(Q)/G_{\bn} \right)  & 
\end{tikzcd}
\]
in which the vertical maps are given by averaging over the Weyl group (i.e. the operation denoted by $\text{Sym}$ in \eqref{eqn:shuffle product}), and so they are surjective. Thus, the fact that any $E \in \text{Im }\jmath_*$ satisfies condition \eqref{eqn:rigorous} follows from the fact that any $E \in \text{Im }\bar{\jmath}_*$ satisfies the same condition. To prove the latter claim, we note that the composition of the two top horizontal maps is $0$ due to \eqref{eqn:choice 2}. Thus, it remains to show that
\begin{equation}
\label{eqn:kappa}
\kappa^*\Big(E(z_{0,1},z_{1,1},z_{1,2},z_{1,3},z_{2,1},z_{2,2},z_{3,1},\dots )\Big) = 
\end{equation} 
\[
\text{proj} \left( E\Big|_{z_{3,1} \mapsto z_{2,1}q^{\mp 1}} \Big|_{z_{1,2} \mapsto z_{2,1}q^{\mp 1}} \Big|_{z_{2,1} \mapsto z_{1,1} q^{\mp 1}} \Big|_{z_{1,1} \mapsto z_{0,1} q^{\pm 2}} \Big|_{z_{0,1} \mapsto z_{1,3} q^{\pm 2}} \right)
\]
for any Laurent polynomial $E$, where $\text{proj}$ denotes projection onto the quotient ring
\begin{equation}
\label{eqn:quotient}
\BK[z_{2,2}^{\pm 1},\dots] \Big/ (z_{2,2} - z_{1,3} q^{\pm 1})^2
\end{equation}
Indeed, asking that the second line of \eqref{eqn:kappa} vanishes is precisely equivalent to condition \eqref{eqn:rigorous}, once we recall the relabeling of our variables $z_{i,u} \leftrightarrow z_{N-i,u}$. To prove formula \eqref{eqn:kappa}, we recall the description of $L$ in \eqref{eqn:choice 1}. Identifying $U$ with the subset $\{x_1=x_2=x_3=x_4=x_5=1\} \subset L$ induces an isomorphism of stacks 
\begin{equation}
\label{eqn:restriction}
L/H \xrightarrow{\sim} U/H'
\end{equation}
where $H' \subset H$ is the subtorus that fixes the coefficients $x_1,\dots,x_5$ in \eqref{eqn:phi}. If we recall that $T$ acts on $\phi_{ij}$ by the integer powers of the character $q$ that appear in the decorated Dynkin diagram \eqref{eqn:dynkin cc}, and that $H$ acts on the $v_{i,u}$'s by the characters $z_{i,u}$, then we see that the operation \eqref{eqn:restriction} acts on $K$-theory by inputing any Laurent polynomial $E$ and outputing the five-fold specialization on the second line of \eqref{eqn:kappa}. Meanwhile, the fact that restricting a $K$-theory class to $U/H'$ is the same as projection onto the quotient \eqref{eqn:quotient} follows from the Claim below (wherein $t = \frac {z_{2,2}}{z_{1,3}q^{\pm 1}}$).

\medskip 

\begin{claim}

Consider $\BA^4$ with coordinates $(a,b,c,d)$, which are acted on by the torus $\BC^*$ with characters $(t,t^{-1},t^{-1},t)$. If we let $Z \hookrightarrow \BA^4$ be the hypersurface
\begin{equation}
\label{eqn:eqn}
\gamma_1 ab + \gamma_2 ac + \gamma_3 bd + \gamma_4 cd = 0
\end{equation}
i.e. a quadric that is smooth away from the origin (we assume $\gamma_1\gamma_4 \neq \gamma_2\gamma_3$) and $U = \BA^4 \backslash Z$ be its open complement, then the exact sequence on Grothendieck groups $G = G_0$ of $\BC^*$-equivariant coherent sheaves
\[
G^{\BC^*}(Z) \xrightarrow{\tau} G^{\BC^*}(\BA^4) \rightarrow G^{\BC^*}(U) \rightarrow 0
\]
has the property that 
\[
\emph{Im }\tau = (1-t)^2 \BZ[t^{\pm 1}] \subset \BZ[t^{\pm 1}] = G^{\BC^*}(\BA^4)
\]
Therefore, restricting a $\BC^*$-equivariant $K$-theory class from $\BA^4$ to $U$ is the same operation as projecting a Laurent polynomial onto the quotient $\BZ[t^{\pm 1}]/(1-t)^2$.

\end{claim}

\medskip 

\noindent Let us prove the Claim: consider the plane $\BA^2 \hookrightarrow \BA^4$ corresponding to the equations $b = c = 0$. Since this plane is contained in $Z$ and its image in $G^{\BC^*}(\BA^4)$ is given by $(1-t)^2$, then we conclude that $(1-t)^2 \in \text{Im }\tau$. To conclude the proof, it suffices to prove that any element $x \in \text{Im }\tau$ is a multiple of $(1-t)^2$. To this end, we consider the commutative diagram with the horizontal maps given by direct image and the vertical maps given by restriction to open subsets:
\[
\begin{tikzcd}
G^{\BC^*}(Z) \arrow[r,"\tau"] \arrow[d,"\rho"'] & G^{\BC^*}(\BA^4) \arrow[d,"\rho'"] \\
G^{\BC^*}(Z \backslash \BA^2) \arrow[r,"\tau'"] & G^{\BC^*}(\BA^4 \backslash \BA^2)
\end{tikzcd}
\]
Since $(Z \backslash \BA^2)/\BC^*$ and $(\BA^4 \backslash \BA^2)/\BC^*$ are easily seen to be isomorphic to the total space of vector bundles on $\BP^1$ (the line bundle $\CO(-2)$ and the vector bundle $\CO(-1) \oplus \CO(-1)$, respectively), the two Grothendieck groups on the bottom of the diagram above are both $\BZ[t^{\pm 1}]/(1-t)^2$. In order to conclude the proof, we need to show that $\rho' \circ \tau = 0$, and for this it suffices to show that $\tau' = 0$. However, since the $G$-theory of total spaces of vector bundles on $\BP^1$ is spanned by the classes $[\CO]$ and $[\CO(1)]$, it suffices to show that $\tau'([\CO_{Z\backslash \BA^2}])=0$. This holds because $\CO_{Z\backslash \BA^2}$ admits the equivariant Koszul resolution
\[
\CO_{\BA^4\backslash \BA^2} \xrightarrow{\text{equation \eqref{eqn:eqn}}} \CO_{\BA^4\backslash \BA^2}
\]
in the ambient space $\BA^4\backslash \BA^2$, and the equation \eqref{eqn:eqn} is $\BC^*$-invariant (thus, the alternating sum of its Koszul complex is $1-1=0$ in $G$-theory). \end{proof}

\end{document}